\documentclass[11pt,a4paper]{article}
\usepackage[T1]{fontenc}
\usepackage{lmodern,microtype}
\usepackage[margin=27mm]{geometry}
\usepackage{amsmath,amssymb,amsthm,mathtools,booktabs}
\usepackage{enumitem,needspace}
\usepackage{xcolor}
\definecolor{linkblue}{RGB}{25,55,105}
\usepackage[colorlinks=true,linkcolor=linkblue,citecolor=linkblue,urlcolor=linkblue]{hyperref}
\hypersetup{
 pdftitle={Integral Exponential Series: Fibonacci Weights and Implicit Equations},
 pdfauthor={Jinping Li, Yueming Zhong, Yangbiao Zhou}
}
\numberwithin{equation}{section}
\newtheorem{theorem}{Theorem}[section]
\newtheorem{lemma}[theorem]{Lemma}
\newtheorem{corollary}[theorem]{Corollary}
\newtheorem{proposition}[theorem]{Proposition}
\theoremstyle{definition}
\newtheorem{example}[theorem]{Example}
\newcommand{\Z}{\mathbb Z}
\newcommand{\Q}{\mathbb Q}
\newcommand{\Fp}{\mathbb F_p}
\newcommand{\EE}{\mathcal E}
\newcommand{\Int}{\operatorname{Int}}
\newcommand{\ord}{\operatorname{ord}}
\title{Integral Exponential Series:\\Fibonacci Weights and Implicit Equations}
\author{%
\begin{tabular}{ccc}
Jinping Li & Yueming Zhong & Yangbiao Zhou\\[0.2em]
\footnotesize\href{mailto:Ljp@hainanu.edu.cn}{\nolinkurl{Ljp@hainanu.edu.cn}} &
\footnotesize\href{mailto:zhongyueming107@gmail.com}{\nolinkurl{zhongyueming107@gmail.com}} &
\footnotesize\href{mailto:yangbiao_zhou@163.com}{\nolinkurl{yangbiao_zhou@163.com}}
\end{tabular}\\[0.75em]
\small Key Laboratory of Engineering Modeling and Statistical Computation of Hainan Province\\
\small School of Mathematics and Statistics, Hainan University, Haikou 570228, P.R. China
}
\date{}
\begin{document}
\maketitle
\begin{abstract}
We apply classical prime-power congruence criteria to two families of
exponential series. For the Fibonacci numbers $F_n$, we characterize the polynomials
$P\in\Z[T]$ for which $\exp(\sum_{n\ge1}P(F_n)x^n/n)$ has integer
coefficients: this holds if and only if $P$ is even and
$P(0)\equiv P(1)\equiv P(2)\pmod5$. In particular, $5$ is the least
positive multiplier for each positive even power of $F_n$, and the series
with weights $5F_n^{2k}$ has degree $5^{k-1}$ over $\Q(x)$ for $k\ge1$.
For $R_n(y)\in\Z[y]$ satisfying $R_{pn}(y)\equiv R_n(y)\pmod p$
for every prime $p$ and every $n\ge1$,
we prove that $H(x,y)=\exp(\sum_{n\ge1}R_n(y)^n x^n/n)$ belongs to
$\Z[y][[x]]$ and that $A(x)=H(x,A(x))$ has a unique integral solution.
For $q\in\Z$, taking $R_n(y)=q^n+y$ proves the conjecture in OEIS A163138;
the Fibonacci result proves the conjecture in A207969 and A207970.
For the implicit family with indeterminate $q$, the coefficient of $x^n$
has degree $n^2$, leading coefficient $1/n$, and nonnegative integer
coefficients in the binomial basis, for every $n\ge1$.
\end{abstract}

\section{Introduction}
All series are formal. Integrality refers to ordinary coefficients in
$\Z$ or in the specified polynomial ring over $\Z$. Congruences between
polynomials or series are coefficientwise. We write $[x^n]F(x)$ for the
coefficient of $x^n$ in $F(x)$.

For a sequence of integer weights $b_n$, integrality of
$\exp(\sum_{n\ge1}b_nx^n/n)$ is governed by congruences comparing
$b_{mp^r}$ and $b_{mp^{r-1}}$. For polynomial weights $b_n(y)$, the
corresponding comparison also involves the substitution $y\mapsto y^p$.
We use these two forms of the classical Dieudonn\'e--Dwork and Dold--Euler
criteria to study Fibonacci weights and an implicit generating function.
These criteria are recalled in Section~\ref{sec:criterion}; see
Byszewski, Graff, and Ward \cite{BGW} and Krattenthaler and M\"uller \cite{KM}.

Let $F_0=0$, $F_1=1$, and $F_{n+2}=F_{n+1}+F_n$. For $P\in\Z[T]$, put
\begin{equation}\label{eq:EP}
 \EE_P(x)=\exp\left(\sum_{n\ge1}P(F_n)\frac{x^n}{n}\right).
\end{equation}
Paul D. Hanna's OEIS entries A207969 and A207970 record the conjecture that
\begin{equation}\label{eq:Gk}
 G_k(x)=\exp\left(\sum_{n\ge1}5F_n^{2k}\frac{x^n}{n}\right)
 \in\Z[[x]]\qquad(k\in\Z_{\ge0});
\end{equation}
the entries display $k=2$ and $k=3$, respectively \cite{OEIS969,OEIS970}.
Our first result determines all such polynomial weights.

\begin{theorem}\label{thm:fibonacci}
For $P\in\Z[T]$, the following statements are equivalent.
\begin{enumerate}[label=\textup{(\roman*)},leftmargin=2.5em]
\item $\EE_P(x)\in\Z[[x]]$.
\item $P(T)=P(-T)$ and $P(0)\equiv P(1)\equiv P(2)\pmod5$.
\item There exist $c\in\Z$ and $U,V\in\Z[Y]$ such that
\begin{equation}\label{eq:polynomial-form}
 P(T)=c+5T^2U(T^2)+(T^6-T^2)V(T^2).
\end{equation}
\end{enumerate}
In \textup{(iii)}, necessarily $c=P(0)$.
\end{theorem}

The second family comes from Hanna's conjecture in OEIS A163138 \cite{OEIS}:
for every integer $q$, the equation
\begin{equation}\label{eq:hanna}
 A(x,q)=\exp\left(\sum_{n\ge1}(q^n+A(x,q))^n\frac{x^n}{n}\right)
\end{equation}
has an integer-coefficient solution. The entry displays the case $q=2$.
We replace $A(x,q)$ in the exponent by an independent variable $y$ and
prove integrality before substitution. The following sufficient condition
follows from the polynomial Dwork criterion and elementary power lifting.

\begin{theorem}\label{thm:implicit}
Let $R_n(y)\in\Z[y]$ for $n\ge1$, and assume that
\begin{equation}\label{eq:weak}
 R_{pn}(y)\equiv R_n(y)\pmod p
 \qquad\text{for every prime $p$ and every $n\ge1$}.
\end{equation}
Then
\begin{equation}\label{eq:H-general}
 H(x,y)=\exp\left(\sum_{n\ge1}R_n(y)^n\frac{x^n}{n}\right)
 \in1+x\Z[y][[x]].
\end{equation}
The equation $A(x)=H(x,A(x))$ has a unique solution in
$1+x\Q[[x]]$, and that solution belongs to $1+x\Z[[x]]$.
\end{theorem}

For fixed $q\in\Z$, Fermat's congruence $q^{pn}\equiv q^n\pmod p$
shows that $R_n(y)=q^n+y$ satisfies \eqref{eq:weak}. Thus
Theorem~\ref{thm:implicit} proves the conjecture for every integer $q$.
The arithmetic arguments in the two applications differ: Fibonacci
weights require separate treatment at $5$, whereas the exponent $n$ in
$R_n(y)^n$ lifts congruences modulo $p$ to the required prime powers.
The implicit equation then follows by coefficientwise iteration.

For every $P\in\Z[T]$, the sequence $P(F_n)$ satisfies a linear recurrence
with constant coefficients. The Fibonacci problem therefore falls within
the setting of Minton's classification of linear recurrence sequences
satisfying Gauss congruences \cite{Minton}. Beukers, Houben, and Straub
\cite[Theorem~8.1]{BHS} give a rational-function formulation and another
proof of that classification, with congruences at all but finitely many
primes. Theorem~\ref{thm:fibonacci} gives an explicit condition on $P$
for integrality at every prime, including the restriction at $5$.
Section~\ref{sec:fibonacci} gives a direct proof. Moss and Ward
\cite{MW} study the related sequence $5F_{n^{2k}}$, obtained by changing
the index, whereas our weights are $5(F_n)^{2k}$. Their use of Fibonacci
congruences at infinitely many primes is relevant to our necessity argument.

Several formulas for the Fibonacci family are already known. Hanna's
formula for $G_1$ appears in A054888 \cite{OEIS888}, and Bala's formulas
for $G_2^5$ and $G_3^{25}$ appear in A207969 and A207970
\cite{OEIS969,OEIS970}. Section~\ref{sec:algebraic} gives their uniform
extension to every $k\ge1$ and proves the exact degree $5^{k-1}$.

Sections~\ref{sec:implicit} and \ref{sec:hanna} treat the implicit
equations. With $q$ an indeterminate, the coefficients of $A(x,q)$
need not lie in $\Z[q]$. Section~\ref{sec:parameter} shows that they
are integer-valued polynomials and determines their degree, leading
coefficient, and nonnegative integer expansion in the binomial basis.

\section{Classical integrality criteria}\label{sec:criterion}
We recall the polynomial form of the Dieudonn\'e--Dwork criterion and
its consequence for exponential series. For integer weights, the
result is the classical Dold--Euler criterion. Proofs are included
to fix the coefficient rings and the Frobenius substitution used below.
For a prime $p$, let
\[
 \Z_{(p)}=\{a/b:a,b\in\Z,\ p\nmid b\}.
\]
Write $v_p(u)$ for the $p$-adic valuation of a nonzero rational number
$u$, normalized by $v_p(p)=1$, and let $\mu$ denote the M\"obius function.
For a finite tuple of indeterminates
$\mathbf y=(y_1,\ldots,y_s)$, let $\sigma_p$ be the $\Q$-algebra
endomorphism of $\Q[\mathbf y]$ defined by $\sigma_p(y_i)=y_i^p$.
On series we use the substitution
$(x,\mathbf y)\mapsto(x^p,\sigma_p(\mathbf y))$, where
$\sigma_p(\mathbf y)=(y_1^p,\ldots,y_s^p)$.

\begin{lemma}[Polynomial form of the Dwork criterion]\label{lem:dwork}
Let $R=\Z_{(p)}[\mathbf y]$ and
$F(x,\mathbf y)\in1+x\Q[\mathbf y][[x]]$. Then
\begin{equation}\label{eq:dwork}
 F\in R[[x]]
 \quad\Longleftrightarrow\quad
 \frac{F(x,\mathbf y)^p}{F(x^p,\sigma_p(\mathbf y))}
      \in1+pxR[[x]].
\end{equation}
\end{lemma}
\begin{proof}
If $F\in R[[x]]$, then
$F(x,\mathbf y)^p\equiv F(x^p,\sigma_p(\mathbf y))\pmod{pR[[x]]}$,
since $\Z_{(p)}/p\Z_{(p)}=\mathbb F_p$. The denominator is a unit
in $R[[x]]$, so division proves the forward implication.

For the converse, write $F=1+\sum_{n\ge1}f_n(\mathbf y)x^n$ and
assume by induction that $f_j\in R$ for $j<N$. Let
$K=1+\sum_{j=1}^{N-1}f_jx^j$. The quotient in \eqref{eq:dwork} for
$U=F/K$ is the quotient for $F$ divided by that for $K$.
The forward implication applies to $K$, so this ratio belongs to
$1+pxR[[x]]$. But
\[
 U=1+f_Nx^N+O(x^{N+1}),\qquad
 U(x^p,\sigma_p(\mathbf y))=1+O(x^{pN}).
\]
The coefficient of $x^N$ in the quotient for $U$ is $pf_N$, which
belongs to $pR$. Hence $f_N\in R$, completing the induction.
\end{proof}

\begin{proposition}\label{prop:ghost}
Suppose $b_n(\mathbf y)\in\Z[\mathbf y]$ and, for every prime $p$,
every $r\ge1$, and every $m\ge1$ with $p\nmid m$,
\begin{equation}\label{eq:ghost}
 b_{mp^r}(\mathbf y)
  \equiv b_{mp^{r-1}}(\sigma_p(\mathbf y))\pmod{p^r}.
\end{equation}
Then
\[
 \exp\left(\sum_{n\ge1}b_n(\mathbf y)\frac{x^n}{n}\right)
 \in1+x\Z[\mathbf y][[x]].
\]
\end{proposition}
\begin{proof}
Fix $p$, set $R=\Z_{(p)}[\mathbf y]$, and write
$S=\sum_{n\ge1}b_n(\mathbf y)x^n/n$. If $p\nmid n$, the coefficient
of $x^n$ in $pS-S(x^p,\sigma_p(\mathbf y))$ is $pb_n/n\in pR$.
For $n=mp^r$ with $p\nmid m$, it is
\[
 \frac{b_{mp^r}(\mathbf y)-b_{mp^{r-1}}(\sigma_p(\mathbf y))}
      {mp^{r-1}}\in pR.
\]
Thus the difference belongs to $pxR[[x]]$. If $V\in xR[[x]]$,
then $\exp(pV)\in1+pxR[[x]]$, because
$v_p(p^j/j!)\ge1$ for every $j\ge1$.
Apply Lemma~\ref{lem:dwork} to $\exp S$ and then intersect the
coefficient rings over all primes, using
$\bigcap_p\Z_{(p)}=\Z$.
\end{proof}

For weights independent of $\mathbf y$, the Frobenius substitution has no effect.
We will also need the converse implication in this case.

\begin{corollary}[Dold--Euler criterion]\label{lem:dold}
Let $b_n\in\Z$ for $n\ge1$. The following conditions are equivalent:
\begin{align}
 &\exp\left(\sum_{n\ge1}b_n\frac{x^n}{n}\right)\in\Z[[x]],
 \label{eq:criterion-int}\\
 &b_{mp^r}\equiv b_{mp^{r-1}}\pmod{p^r}
 \quad\text{for every prime $p$, $r\ge1$, and $m\ge1$ with $p\nmid m$}.
 \label{eq:criterion-cong}
\end{align}
\end{corollary}
\begin{proof}
For sufficiency, take the constant polynomials $b_n(y)=b_n$ in
Proposition~\ref{prop:ghost}.

Conversely, an integral series with constant term $1$ has a unique Euler
product with integer exponents. At stage $n$, if the remaining series
is $1+c_nx^n+O(x^{n+1})$, set $e_n=c_n$ and multiply by
$(1-x^n)^{e_n}$. Taking logarithms gives
\begin{equation}\label{eq:euler}
 \exp\left(\sum_{n\ge1}b_n\frac{x^n}{n}\right)
 =\prod_{n\ge1}(1-x^n)^{-e_n},
 \qquad b_n=\sum_{d\mid n}de_d.
\end{equation}
Equivalently, by M\"obius inversion,
\begin{equation}\label{eq:euler-exponent}
 e_n=\frac1n\sum_{d\mid n}\mu(n/d)b_d.
\end{equation}
For $p\nmid m$ it follows that
\[
 b_{mp^r}-b_{mp^{r-1}}=p^r\sum_{d\mid m}de_{p^rd},
\]
which proves \eqref{eq:criterion-cong}.
\end{proof}

\section{Polynomial Fibonacci weights}\label{sec:fibonacci}
\subsection{Fibonacci congruences}
Set
\[
 \alpha=\frac{1+\sqrt5}{2},\qquad
 \beta=\frac{1-\sqrt5}{2},\qquad L_n=\alpha^n+\beta^n\quad(n\ge0).
\]
Thus $F_n=(\alpha^n-\beta^n)/(\alpha-\beta)$, $\alpha\beta=-1$,
and $L_n\in\Z$.

\begin{lemma}\label{lem:squares}
For every prime $p\ne5$, every $r\ge1$, and every $m\ge1$ with $p\nmid m$,
\begin{equation}\label{eq:square-cong}
 F_{mp^r}^2\equiv F_{mp^{r-1}}^2\pmod{p^r}.
\end{equation}
Also, $5^s\mid F_{m5^s}$ for all $m\ge1$ and $s\ge0$.
\end{lemma}
\begin{proof}
Binet's formula gives $5F_n^2=L_{2n}-2(-1)^n$. This gives the
identity recorded by Hanna in A054888 \cite{OEIS888}:
\begin{equation}\label{eq:base-rational}
 \exp\left(\sum_{n\ge1}5F_n^2\frac{x^n}{n}\right)
   =\frac{(1+x)^2}{(1-\alpha^2x)(1-\beta^2x)}
   =\frac{(1+x)^2}{1-3x+x^2}\in\Z[[x]].
\end{equation}
Corollary~\ref{lem:dold} gives the prime-power congruences for $5F_n^2$;
cancelling $5$ proves \eqref{eq:square-cong}.

Factoring $\alpha^{5n}-\beta^{5n}$ and using $\alpha\beta=-1$ gives
\[
 F_{5n}=F_n\bigl(L_n^4-3(-1)^nL_n^2+1\bigr).
\]
Substituting $L_n^2=5F_n^2+4(-1)^n$ yields
\begin{equation}\label{eq:fivefold}
 F_{5n}=5F_n\bigl(5F_n^4+5(-1)^nF_n^2+1\bigr).
\end{equation}
Iteration proves $5^s\mid F_{m5^s}$. The factor in parentheses is
$1$ modulo $5$, so $v_5(F_{5n})=v_5(F_n)+1$ for $n\ge1$.
\end{proof}

\begin{lemma}\label{lem:inert}
There are infinitely many odd primes $p$ for which $5$ is a quadratic
nonresidue modulo $p$. At each such prime,
\begin{equation}\label{eq:inert-cong}
 F_{mp}\equiv-F_m\pmod p\qquad(m\ge1).
\end{equation}
\end{lemma}
\begin{proof}
With $(a/p)$ denoting the Legendre symbol, quadratic reciprocity gives
$(5/p)=(p/5)$ for odd $p\ne5$.
To see infinitude, let $S$ be any finite set of odd primes different
from $5$. By the Chinese remainder theorem choose a positive integer
$N>1$ with $N\equiv2\pmod5$ and
$N\equiv1\pmod{2\prod_{\ell\in S}\ell}$. The product of the symbols
$(\ell/5)$ over the prime factors of $N$, counted with multiplicity,
is $(N/5)=-1$. Thus some odd prime factor $p$ of $N$ satisfies
$(p/5)=-1$, and this prime does not lie in $S$.

For such a prime, $z^2-z-1$ is irreducible over $\Fp$. Let $a,b$ be its
roots in $\mathbb F_{p^2}$. Frobenius exchanges them, so $a^p=b$ and
$b^p=a$. Since $a-b\ne0$, the Fibonacci recurrence in this field yields
\[
 F_{mp}=\frac{a^{mp}-b^{mp}}{a-b}
       =\frac{b^m-a^m}{a-b}=-F_m.
\]
Both sides belong to $\Fp$, proving the congruence.
\end{proof}

\subsection{Proof of the classification}
\begin{proof}[Proof of Theorem~\ref{thm:fibonacci}]
Suppose first that $\EE_P$ is integral. For every fixed $m\ge1$,
Corollary~\ref{lem:dold} and Lemma~\ref{lem:inert} imply
\[
 P(-F_m)\equiv P(F_m)\pmod p
\]
for infinitely many primes $p$; the finitely many divisors of $m$ may be
excluded. Hence the fixed integer $P(-F_m)-P(F_m)$ is zero. The numbers
$F_m$ assume infinitely many distinct values, so the polynomial
$P(-T)-P(T)$ vanishes identically. Write $P(T)=Q(T^2)$ with $Q\in\Z[Y]$.

At $p=5$, the cases $m=1$ and $m=3$ of Corollary~\ref{lem:dold} give
\[
 P(F_5)\equiv P(F_1),\qquad P(F_{15})\equiv P(F_3)\pmod5.
\]
As $F_1=1$, $F_3=2$, and $5\mid F_5,F_{15}$, these are precisely
$P(0)\equiv P(1)\equiv P(2)\pmod5$. Thus (i) implies (ii).

Assume (ii), and again write $P(T)=Q(T^2)$. For $p\ne5$,
Lemma~\ref{lem:squares} and evaluation of $Q$ give all congruences
\eqref{eq:criterion-cong} for $b_n=P(F_n)$.
For $p=5$ and $r=1$, every square modulo $5$ is $0$, $1$, or $4$,
and (ii) says that $Q$ has the same value at these three residues.
Consequently
\[
 P(F_{5m})\equiv Q(0)\equiv P(F_m)\pmod5.
\]
For $r\ge2$, Lemma~\ref{lem:squares} gives
\[
 5^{2(r-1)}\mid F_{m5^{r-1}}^2,
 \qquad 5^{2r}\mid F_{m5^r}^2.
\]
Since $2(r-1)\ge r$, both $Q(F_{m5^{r-1}}^2)$ and
$Q(F_{m5^r}^2)$ are congruent to $Q(0)$ modulo $5^r$.
Corollary~\ref{lem:dold} now proves (i).

Finally, in $\mathbb F_5[Y]$, condition (ii) is equivalent to
\[
 Y(Y-1)(Y+1)=Y^3-Y\quad\text{dividing}\quad Q(Y)-Q(0).
\]
Lift the quotient to a polynomial $V\in\Z[Y]$. The difference
$Q(Y)-Q(0)-(Y^3-Y)V(Y)$ is divisible by $5$ coefficientwise and
vanishes at $Y=0$, so it equals $5YU(Y)$ for some $U\in\Z[Y]$.
Substitution $Y=T^2$ gives (iii). The reverse implication is immediate
from the same three residues.
\end{proof}

\begin{corollary}\label{cor:scalar}
Let $c\in\Z$ and $j\in\Z_{\ge0}$. Then
\[
 \exp\left(\sum_{n\ge1}cF_n^j\frac{x^n}{n}\right)\in\Z[[x]]
\]
holds under exactly the following conditions:
\begin{enumerate}[label=\textup{(\alph*)},leftmargin=2.5em]
\item $j=0$: every $c\in\Z$ is allowed;
\item $j>0$ is even: $5\mid c$;
\item $j$ is odd: $c=0$.
\end{enumerate}
In particular, \eqref{eq:Gk} holds for every $k\ge0$, and $5$ is the
least positive multiplier for every positive even power.
\end{corollary}
\begin{proof}
Apply Theorem~\ref{thm:fibonacci} to $P(T)=cT^j$. For $j=0$ the series is
$(1-x)^{-c}$. For positive even $j$, comparison of $P(0)$ and $P(1)$
requires $5\mid c$, which is also sufficient. A nonzero odd monomial
is not an even polynomial.
\end{proof}

\begin{example}
The polynomial $T^6-T^2$ satisfies the criterion without an overall
factor $5$, and therefore
\[
 \exp\left(\sum_{n\ge1}(F_n^6-F_n^2)\frac{x^n}{n}\right)\in\Z[[x]].
\]
In contrast, $P(T)=T^4-T^2$ has $P(0)=P(1)=0$ but
$P(2)=12\not\equiv0\pmod5$, so its exponential is not integral.
The scalar condition in Corollary~\ref{cor:scalar} therefore cannot be
applied separately to each monomial of an arbitrary polynomial weight.
\end{example}

\section{Rational powers and exact algebraic degree}\label{sec:algebraic}
Recall that $L_n=\alpha^n+\beta^n$ is the Lucas number.

\begin{theorem}\label{thm:degree}
For an integer $k\ge1$, let $d_k=5^{k-1}$ and define
\begin{equation}\label{eq:Rk}
\begin{split}
 \mathcal R_k(x)={}&\bigl(1-(-1)^kx\bigr)^{(-1)^{k+1}\binom{2k}{k}}\\
 &\times\prod_{j=0}^{k-1}
 \bigl(1-(-1)^jL_{2(k-j)}x+x^2\bigr)^%
 {(-1)^{j+1}\binom{2k}{j}}.
\end{split}
\end{equation}
Then $\mathcal R_k\in\Q(x)$, $G_k^{d_k}=\mathcal R_k$, and
\begin{equation}\label{eq:degree}
 [\Q(x,G_k):\Q(x)]=5^{k-1}.
\end{equation}
Consequently $G_1$ is rational and every $G_k$ with $k\ge2$ is
algebraic and irrational over $\Q(x)$.
\end{theorem}
\begin{proof}
The binomial theorem gives
\[
 5^kF_n^{2k}=
 \sum_{j=0}^{2k}(-1)^j\binom{2k}{j}
       \bigl(\alpha^{2k-j}\beta^j\bigr)^n.
\]
Taking formal exponentials shows that $G_k^{d_k}$ is the product of
$(1-\alpha^{2k-j}\beta^jx)^{-(-1)^j\binom{2k}{j}}$ over
$0\le j\le2k$. Pairing $j$ with $2k-j$ gives the quadratic factors
in \eqref{eq:Rk}; the unpaired term $j=k$ gives its linear factor.

Put $P_k(x)=1-L_{2k}x+x^2$. Its discriminant is
\[
 L_{2k}^2-4=5F_{2k}^2,
\]
which is not a rational square. Thus $P_k$ is irreducible over $\Q$.
The linear factors before pairing have parameters
$\lambda_j=(-1)^j\alpha^{2k-2j}$. Since $\alpha>1$, their absolute
values, and hence the parameters themselves, are distinct.
It follows that the factor $P_k$ in $\mathcal R_k$ has exponent exactly $-1$,
and none of the other factors contributes to its valuation.
Writing $\ord_{P_k}$ for the exponent of $P_k$ in a rational function,
we obtain
\[
 \ord_{P_k}(\mathcal R_k^{-1})=1.
\]
For $k=1$, the product formula gives $G_1\in\Q(x)$. If $k\ge2$, the
polynomial $Y^{d_k}-\mathcal R_k^{-1}$ is Eisenstein in $\Q[x]_{(P_k)}$, the
localization at $(P_k)$: its leading coefficient is $1$, its intermediate
coefficients are zero, and its constant coefficient has valuation $1$.
It is therefore irreducible over $\Q(x)$. Since $G_k^{-1}$ is a root,
\eqref{eq:degree} follows.
\end{proof}

For example,
\begin{align*}
 G_1(x)&=\frac{(1+x)^2}{1-3x+x^2},\\
 G_2(x)^5&=\frac{(1+3x+x^2)^4}{(1-7x+x^2)(1-x)^6},\\
 G_3(x)^{25}&=
 \frac{(1+x)^{20}(1+7x+x^2)^6}
 {(1-3x+x^2)^{15}(1-18x+x^2)}.
\end{align*}
The first identity is recorded in A054888 \cite{OEIS888}; the other two
are Bala's formulas \cite{OEIS969,OEIS970}. The simple poles in the
last two formulas give the degrees $5$ and $25$ by the same Eisenstein argument.

\section{Integral implicit equations}\label{sec:implicit}
\begin{lemma}[Power lifting]\label{lem:lift}
Let $p$ be a prime. If $U,V\in\Z[\mathbf y]$ and $U\equiv V\pmod p$, then
\begin{equation}\label{eq:lift}
 U^{mp^{r-1}}\equiv V^{mp^{r-1}}\pmod{p^r}
 \qquad(m,r\in\Z_{\ge1}).
\end{equation}
\end{lemma}
\begin{proof}
For $t\ge1$, the binomial theorem shows that
$U\equiv V\pmod{p^t}$ implies $U^p\equiv V^p\pmod{p^{t+1}}$.
Indeed, write $U=V+p^tW$ with $W\in\Z[\mathbf y]$.
Every term of $(V+p^tW)^p-V^p$ is divisible by $p^{t+1}$:
the intermediate binomial coefficients are divisible by $p$, and
$tp\ge t+1$. Iterating and then taking an $m$-th power proves
\eqref{eq:lift}, including when $p=2$.
\end{proof}

\begin{proof}[Proof of Theorem~\ref{thm:implicit}]
Set $b_n(y)=R_n(y)^n$. Fix $p$ and write $n=mp^r$, where $p\nmid m$
and $r\ge1$. By \eqref{eq:weak} and Frobenius modulo $p$,
\[
 R_n(y)^p\equiv R_{n/p}(y)^p\equiv R_{n/p}(y^p)\pmod p.
\]
Raising these polynomials to the power $n/p=mp^{r-1}$ and applying
Lemma~\ref{lem:lift} gives
\[
 b_n(y)=\bigl(R_n(y)^p\bigr)^{n/p}
       \equiv R_{n/p}(y^p)^{n/p}=b_{n/p}(y^p)\pmod{p^r}.
\]
Proposition~\ref{prop:ghost} proves \eqref{eq:H-general}.

Write $H(x,y)=1+\sum_{n\ge1}h_n(y)x^n$ with $h_n\in\Z[y]$.
Each coefficient is a polynomial, so substitution of an arbitrary
$B(x)\in\Q[[x]]$ for $y$ is meaningful even when $B(0)\ne0$.
Since each $h_n(B)-h_n(C)$ is divisible by $B-C$, for
$B,C\in\Q[[x]]$ and $s\ge0$ we have
\begin{equation}\label{eq:contraction}
 B\equiv C\pmod{x^s}
 \quad\Longrightarrow\quad
 H(x,B)\equiv H(x,C)\pmod{x^{s+1}}.
\end{equation}
Starting with $A_0=1$, define $A_{j+1}=H(x,A_j)$. All iterates lie
in $1+x\Z[[x]]$, and \eqref{eq:contraction} yields
$A_{j+1}-A_j\in x^{j+1}\Z[[x]]$. Their coefficients therefore
stabilize to an integral series $A$. Passing to the stabilized coefficients
in $A_{j+1}=H(x,A_j)$ gives $A=H(x,A)$.

Every solution has constant term $1$. If two solutions in
$1+x\Q[[x]]$ first differed at degree $s$, \eqref{eq:contraction} would force their
difference to vanish through degree $s$, a contradiction. This proves
uniqueness in $1+x\Q[[x]]$.
\end{proof}

\begin{corollary}\label{cor:multivariate}
If $R_n(\mathbf y)\in\Z[\mathbf y]$ satisfies
$R_{pn}(\mathbf y)\equiv R_n(\mathbf y)\pmod p$ for every $p,n$, then
\[
 \exp\left(\sum_{n\ge1}R_n(\mathbf y)^n\frac{x^n}{n}\right)
 \in1+x\Z[\mathbf y][[x]].
\]
\end{corollary}
\begin{proof}
The same power-lifting argument gives \eqref{eq:ghost} with
$\sigma_p(y_i)=y_i^p$ for each variable. Apply Proposition~\ref{prop:ghost}.
\end{proof}

\begin{corollary}\label{cor:families}
Let $q_1,\ldots,q_s$ be fixed integers and let
$P(y,z_1,\ldots,z_s)\in\Z[y,z_1,\ldots,z_s]$. Then
\[
 A(x)=\exp\left(\sum_{n\ge1}
       P(A(x),q_1^n,\ldots,q_s^n)^n\frac{x^n}{n}\right)
\]
has a unique solution in $1+x\Q[[x]]$, and this solution is integral.
\end{corollary}
\begin{proof}
Take $R_n(y)=P(y,q_1^n,\ldots,q_s^n)$.
The congruences $q_i^{pn}\equiv q_i^n\pmod p$ verify
\eqref{eq:weak}, so Theorem~\ref{thm:implicit} applies.
\end{proof}

For instance, the result applies to
$R_n(y)=q^n+y^d$ with fixed integers $q\in\Z$ and $d\ge0$, and to
$R_n(y)=u^n+v^n+y-y^2$ with arbitrary integers $u,v$.
No sign condition on the coefficients of $P$ is required.

\section{The equation in A163138}\label{sec:hanna}
\begin{corollary}\label{cor:OEIS}
For every fixed $q\in\Z$,
\begin{equation}\label{eq:Hq}
 H_q(x,y)=\exp\left(\sum_{n\ge1}(q^n+y)^n\frac{x^n}{n}\right)
 \in1+x\Z[y][[x]],
\end{equation}
and \eqref{eq:hanna} has a unique solution in $1+x\Z[[x]]$.
\end{corollary}
\begin{proof}
Apply Theorem~\ref{thm:implicit} with $R_n(y)=q^n+y$.
\end{proof}

The specialization
\[
 H_q(x,0)=\exp\left(\sum_{n\ge1}q^{n^2}\frac{x^n}{n}\right)
\]
is already known to be integral for every integer $q$, by Bala's proof
in A155200 \cite{OEIS5200}. For nonnegative integers $q,m$,
Fried \cite[Theorem~43]{Fried} proves integrality of
$\exp(\sum_{n\ge1}q^{n^m}x^n/n)$.
Corollary~\ref{cor:OEIS} establishes integrality in $\Z[y][[x]]$
before the implicit substitution $y=A(x,q)$.

The first two nonconstant coefficients of $H_q$ are
\begin{equation}\label{eq:H-first}
 H_q(x,y)=1+(q+y)x+
 \left(y^2+(q+q^2)y+\frac{q^2+q^4}{2}\right)x^2+O(x^3).
\end{equation}
For each integer $q$, the expression $(q^2+q^4)/2$ is an integer.

\begin{proposition}\label{prop:rearrange}
With $q$ an indeterminate and $H_q$ defined by the same exponential
formula as in \eqref{eq:Hq}, the following identity holds in
$\Q[q,y][[x]]$:
\begin{equation}\label{eq:rearrange}
 H_q(x,y)=\frac1{1-xy}
 \exp\left(\sum_{r\ge1}
    \frac{q^{r^2}x^r}{r(1-q^rxy)^r}\right).
\end{equation}
Consequently substitution $y=A(x,q)$ gives the alternative equation
recorded in \cite{OEIS}.
\end{proposition}
\begin{proof}
Expand $(q^n+y)^n$ and write $n=r+j$, where $r$ counts the factors
$q^n$ and $j$ counts the factors $y$. The terms with $r=0$ sum to
$-\log(1-xy)$. For $r\ge1$, use
$\binom{r+j}{j}/(r+j)=\binom{r+j-1}{j}/r$ to obtain
\[
 \sum_{j\ge0}\frac{\binom{r+j}{j}}{r+j}
 q^{r(r+j)}x^{r+j}y^j
 =\frac{q^{r^2}x^r}{r}
   \sum_{j\ge0}\binom{r+j-1}{j}(q^rxy)^j
 =\frac{q^{r^2}x^r}{r(1-q^rxy)^r}.
\]
Every coefficient uses finitely many terms. Exponentiating proves
\eqref{eq:rearrange}.
\end{proof}

At $q=0$, equation~\eqref{eq:hanna} is
$A=(1-xA)^{-1}$, equivalently $A=1+xA^2$, so its coefficients are
the Catalan numbers. At $q=1$, it is
$A=(1-x(1+A))^{-1}$, equivalently $A=1+xA+xA^2$, giving the large
Schr\"oder numbers, as already noted in \cite{OEIS}.

\section{Integer-valued parameter polynomials}\label{sec:parameter}
We now regard $q$ as an indeterminate. The coefficientwise contraction
argument constructs a unique solution
\[
 A(x,q)=\sum_{n\ge0}a_n(q)x^n\in1+x\Q[q][[x]]
\]
of \eqref{eq:hanna}. For each $q_0\in\Z$, specialization at $q=q_0$
commutes with the construction and gives the solution of
Corollary~\ref{cor:OEIS} with parameter $q_0$.
Recall that $\Int(\Z)$ denotes the ring of polynomials in $\Q[q]$
that take integer values at every integer. Its binomial-basis description
is classical; see Cahen and Chabert \cite{CC}.

\Needspace{10\baselineskip}
\begin{theorem}\label{thm:parameter}
For every $n\ge1$, the coefficient $a_n(q)$ has the following properties:
\begin{enumerate}[label=\textup{(\roman*)},leftmargin=2.5em]
\item $a_n(q)\in\Int(\Z)$ and all its monomial coefficients are
nonnegative rational numbers;
\item $\deg a_n=n^2$ and $[q^{n^2}]a_n(q)=1/n$;
\item there are nonnegative integers $d_{n,j}$ such that
\begin{equation}\label{eq:newton}
 a_n(q)=\sum_{j=0}^{n^2}d_{n,j}\binom qj,
 \qquad d_{n,n^2}=\frac{(n^2)!}{n}.
\end{equation}
\end{enumerate}
In particular, $a_n(q)\notin\Z[q]$ for every $n\ge2$.
\end{theorem}
\begin{proof}
The integer-valued assertion follows from Corollary~\ref{cor:OEIS}
and compatibility with specialization. To prove the remaining claims,
write
\[
 S(x,q)=\sum_{m\ge1}(q^m+A(x,q))^m\frac{x^m}{m}
       =\sum_{n\ge1}s_n(q)x^n,
 \qquad A=\exp S.
\]
The coefficient $s_n$ only involves $a_i$ with $i<n$, since every
summand in $S$ contains a positive power of $x$. More explicitly,
\begin{equation}\label{eq:sn}
 s_n(q)=\sum_{m=1}^n\frac1m\sum_{j=0}^m
 \binom mj q^{m(m-j)}[x^{n-m}]A(x,q)^j.
\end{equation}
Starting with $a_0=1$, both \eqref{eq:sn} and the expansion of $\exp S$
use nonnegative rational coefficients. Induction proves the nonnegativity
assertion in (i).

For the degree calculation, assume inductively that
$\deg a_i,\deg s_i\le i^2$ for $1\le i<n$, and recall that $a_0=1$.
A product contributing to
$[x^{n-m}]A^j$ has total index $n-m$, so its degree in $q$ is at most
$(n-m)^2$. For $m<n$, every corresponding term of \eqref{eq:sn}
therefore has degree at most
\[
 m(m-j)+(n-m)^2\le m^2+(n-m)^2<n^2.
\]
For $m=n$, the term $j=0$ is $q^{n^2}/n$, and all terms with $j>0$
have smaller degree. Consequently
\[
 s_n(q)=\frac{q^{n^2}}n+\text{a polynomial of degree less than $n^2$}.
\]
In the expansion $A=1+S+S^2/2!+\cdots$, a contribution to $x^n$ from
$S^r$ with $r\ge2$ uses positive indices $i_1+\cdots+i_r=n$.
Its degree is at most $i_1^2+\cdots+i_r^2<n^2$. This proves (ii)
by simultaneous induction for $a_n$ and $s_n$, starting with
$a_1=s_1=q+1$.

Let $\Delta f(q)=f(q+1)-f(q)$. Newton's finite-difference formula gives
\eqref{eq:newton} with $d_{n,j}=\Delta^ja_n(0)$. These are integers
because all the values $a_n(0),a_n(1),\ldots$ are integers. If
$a_n(q)=\sum_e c_eq^e$ with $c_e\ge0$, then
\[
 d_{n,j}=\sum_e c_e\,j!\,\left\{\begin{matrix}e\\j\end{matrix}\right\}\ge0,
\]
where the braces denote Stirling numbers of the second kind.
Finally, the leading coefficient of $\binom q{n^2}$ is $1/(n^2)!$;
part (ii) gives $d_{n,n^2}=(n^2)!/n$. For $n\ge2$, the leading
coefficient $1/n$ also shows that $a_n\notin\Z[q]$.
\end{proof}

The first parameter polynomials are
\begin{align*}
 a_0(q)&=1,\qquad a_1(q)=q+1,\\
 a_2(q)&=\frac{q^4+3q^2+4q+4}{2},\\
 a_3(q)&=\frac{2q^9+6q^6+3q^5+6q^4+19q^3+30q^2+36q+30}{6}.
\end{align*}
For example, the binomial expansion is
\[
 a_2(q)=2+4\binom q1+10\binom q2+18\binom q3+12\binom q4.
\]
It exhibits integer values at negative integers as well, since
$\binom qj\in\Z$ for every integer $q$ and every $j\ge0$.

\section{Coefficient calculations and verification}
\subsection{Fibonacci weights}
If $\EE_P(x)=\sum_{n\ge0}a_P(n)x^n$, formal differentiation of
\eqref{eq:EP} gives the exact recurrence
\begin{equation}\label{eq:recurrence}
 a_P(0)=1,\qquad
 na_P(n)=\sum_{j=1}^nP(F_j)a_P(n-j)\quad(n\ge1).
\end{equation}
This recurrence is convenient for verification, although its division
by $n$ does not itself prove integrality. For the first values of $k$,
the coefficients $[x^n]G_k$ are as follows.
\begin{center}
\begin{tabular}{c r r r r r r}
\toprule
$k\backslash n$&$0$&$1$&$2$&$3$&$4$&$5$\\
\midrule
$0$&1&5&15&35&70&126\\
$1$&1&5&15&40&105&275\\
$2$&1&5&15&60&295&1625\\
$3$&1&5&15&140&1505&21875\\
\bottomrule
\end{tabular}
\end{center}
\subsection{Implicit coefficients}
For fixed $q\in\Z$, write $H_q(x,y)=\sum_{n\ge0}h_n(y)x^n$. Its logarithmic
derivative gives
\begin{equation}\label{eq:h-recur}
 h_0(y)=1,\qquad
 nh_n(y)=\sum_{j=1}^n(q^j+y)^j h_{n-j}(y).
\end{equation}
After these polynomials have been calculated, the coefficients of the
implicit solution satisfy
\begin{equation}\label{eq:a-recur}
 a_0(q)=1,\qquad
 a_n(q)=\sum_{j=1}^n[x^{n-j}]
 h_j\left(\sum_{i=0}^{n-1}a_i(q)x^i\right).
\end{equation}
Equation~\eqref{eq:a-recur} uses only previously determined coefficients.
The first values are shown below.
\begin{center}
\begin{tabular}{c r r r r r r}
\toprule
$q\backslash n$&$0$&$1$&$2$&$3$&$4$&$5$\\
\midrule
$-2$&1&$-1$&12&$-119$&20952&$-5731946$\\
$-1$&1&0&2&2&12&20\\
$0$&1&1&2&5&14&42\\
$1$&1&2&6&22&90&394\\
$2$&1&3&20&329&22584&7938470\\
\bottomrule
\end{tabular}
\end{center}

\subsection{Exact verification}
The accompanying script uses exact rational and integer arithmetic.
For Fibonacci weights, it checks \eqref{eq:recurrence} for
$0\le k\le12$ through degree $100$, verifies \eqref{eq:fivefold} for
$1\le n\le100$, and compares the two sides of \eqref{eq:Rk} for
$1\le k\le5$ through degree $16$. It also checks the integral weights
$T^6-T^2$, $T^8-T^4$, and several combinations from
\eqref{eq:polynomial-form}, together with nonintegral examples.

For the implicit family, the script checks \eqref{eq:h-recur} and
\eqref{eq:a-recur} for $-6\le q\le6$ through degree $14$, and independently
substitutes the resulting series into \eqref{eq:hanna}. Symbolic
calculations check Theorem~\ref{thm:parameter} for $1\le n\le5$,
including all binomial-basis coefficients. A further $217$ checks
test the coefficientwise prime-power congruences, including $p=2$
and negative parameters. These finite checks supplement the proofs.


\begin{thebibliography}{99}
\bibitem{BHS}
F. Beukers, M. Houben, and A. Straub,
Gauss congruences for rational functions in several variables,
\emph{Acta Arith.} \textbf{184} (2018), 341--362.
\href{https://doi.org/10.4064/aa170614-13-7}{doi:10.4064/aa170614-13-7}.

\bibitem{BGW}
J. Byszewski, G. Graff, and T. Ward,
Dold sequences, periodic points, and dynamics,
\emph{Bull. Lond. Math. Soc.} \textbf{53} (2021), 1263--1298.
\href{https://doi.org/10.1112/blms.12531}{doi:10.1112/blms.12531}.

\bibitem{CC}
P.-J. Cahen and J.-L. Chabert,
\emph{Integer-Valued Polynomials}, Mathematical Surveys and Monographs,
vol. 48, American Mathematical Society, Providence, RI, 1997.
\href{https://doi.org/10.1090/surv/048}{doi:10.1090/surv/048}.

\bibitem{Fried}
S. Fried,
Further proofs of conjectures from the OEIS,
preprint, 2026.
\href{https://arxiv.org/abs/2607.24832}{arXiv:2607.24832}.

\bibitem{KM}
C. Krattenthaler and T. W. M\"uller,
Truncated versions of Dwork's lemma for exponentials of power series
and $p$-divisibility of arithmetic functions,
\emph{Adv. Math.} \textbf{283} (2015), 489--529.
\href{https://doi.org/10.1016/j.aim.2015.05.018}{doi:10.1016/j.aim.2015.05.018}.

\bibitem{Minton}
G. T. Minton,
Linear recurrence sequences satisfying congruence conditions,
\emph{Proc. Amer. Math. Soc.} \textbf{142} (2014), 2337--2352.
\href{https://doi.org/10.1090/S0002-9939-2014-12168-X}{doi:10.1090/S0002-9939-2014-12168-X}.

\bibitem{MW}
P. Moss and T. Ward,
Fibonacci along even powers is (almost) realizable,
\emph{Fibonacci Quart.} \textbf{60} (2022), 40--47.
\href{https://arxiv.org/abs/2011.13068}{arXiv:2011.13068}.

\bibitem{OEIS888}
OEIS Foundation Inc.,
\emph{The On-Line Encyclopedia of Integer Sequences}, entry A054888,
contributed by P. Dominici, May 23, 2000;
formula involving Fibonacci squares contributed by P. D. Hanna,
February 21, 2012.
\url{https://oeis.org/A054888}.

\bibitem{OEIS5200}
OEIS Foundation Inc.,
\emph{The On-Line Encyclopedia of Integer Sequences}, entry A155200,
contributed by P. D. Hanna, February 4, 2009;
integrality proof contributed by P. Bala, December 25, 2019.
\url{https://oeis.org/A155200}. 

\bibitem{OEIS}
OEIS Foundation Inc.,
\emph{The On-Line Encyclopedia of Integer Sequences}, entry A163138,
contributed by P. D. Hanna, August 7, 2009.
\url{https://oeis.org/A163138}.

\bibitem{OEIS969}
OEIS Foundation Inc.,
\emph{The On-Line Encyclopedia of Integer Sequences}, entry A207969,
contributed by P. D. Hanna, February 22, 2012;
algebraic formula contributed by P. Bala, April 3, 2014.
\url{https://oeis.org/A207969}.

\bibitem{OEIS970}
OEIS Foundation Inc.,
\emph{The On-Line Encyclopedia of Integer Sequences}, entry A207970,
contributed by P. D. Hanna, February 22, 2012;
algebraic formula contributed by P. Bala, April 3, 2014.
\url{https://oeis.org/A207970}. 

\end{thebibliography}
\end{document}